\documentclass[twoside,12pt]{article}
\usepackage{amssymb,amsmath,bm,euscript}
\usepackage{graphicx,color}
\usepackage{ifpdf}

\usepackage[colorlinks=true]{hyperref}

\newtheorem{proposition}{Proposition}[section]
\newtheorem{theorem}[proposition]{Theorem}

\newtheorem{corollary}[proposition]{Corollary}
\newtheorem{definition}[proposition]{Definition}

\def\S{\mathcal{S}}

\def\R{\mathbb{R}}

\def\I{\mathcal{I}}

\def\Q{\mathcal{Q}}
\def\A{{\mathcal{A}}}
\def\B{\mathcal{B}}
\def\C{\mathcal{C}}
\def\CC{\mathbb{C}}

\def\D{\mathcal{D}}
\def\U{\mathcal{U}}
\def\V{\mathcal{V}}
\def\R{\mathbb{R}}
\def\S{\mathcal{S}}

\def\Y{{\mathcal{Y}}}
\def\Z{\mathcal{Z}}

\def\E{\mathcal{E}}

\def\0{{\bf 0}}

\title{\bf An Orthogonal Equivalence Theorem for Third Order Tensors\footnote{To appear in: Journal of Industrial and Management Optimization.}}
\author{ \hspace{1mm}
	Liqun Qi\thanks{Department of Mathematics, School of Science, Hangzhou Dianzi University, Hangzhou 310018 China;
		({\tt liqun.qi@polyu.edu.hk}). },
Chen Ling\thanks{Department of Mathematics,  Hangzhou Dianzi University, Hangzhou 310018, China; ({\tt macling@hdu.edu.cn}).
 This author's work was supported by Natural Science Foundation of China (No. 11971138) and Natural Science Foundation of Zhejiang Province (No. LY19A010019, LD19A010002).},
Jinjie Liu \thanks{School of Mathematical Sciences, Shanghai Jiao Tong University, Shanghai, 200240, China;({\tt jinjie.liu@sjtu.edu.cn}).
This author's work was supported by Natural Science Foundation of China (No. 12001366, No. 11801479).Corresponding Author. },
\ and\
Chen Ouyang\thanks{School of Computer Science and Technology, Dongguan University of Technology, Dongguan, 523000, China;({\tt
oych26@163.com}). This author's work was supported by Natural Science Foundation of China (No.11971106) and  Guangdong Universities' Special Projects in Key Fields of Natural Science(No. 2019KZDZX1005).}
}

\begin{document}
\date{\today}
\maketitle

\begin{abstract}
In 2011, Kilmer and Martin proposed tensor singular value decomposition (T-SVD) for third order tensors.  Since then, T-SVD has applications in low rank tensor approximation, tensor recovery, multi-view clustering, multi-view feature extraction, tensor sketching, etc.   By going through the Discrete Fourier Transform (DFT), matrix SVD and inverse DFT, a third order tensor is mapped to an f-diagonal third order tensor.   We call this a Kilmer-Martin mapping.  We show that the Kilmer-Martin mapping of a third order tensor is invariant if that third order tensor is taking T-product with some orthogonal tensors.
We define singular values and T-rank of that third order tensor based upon its Kilmer-Martin mapping.   Thus, tensor tubal rank, T-rank, singular values and T-singular values of a third order tensor are invariant when it is taking T-product with some orthogonal tensors.  Some properties of singular values, T-rank and best T-rank one approximation are discussed.

 \noindent {\bf Key words.}  Third order tensors, orthogonal equivalence, singular value, T-rank, the best T-rank one approximation.

\noindent {\bf AMS subject classifications. }{15A69, 15A18}
\end{abstract}

\section{Introduction}
The T-product operation, T-SVD factorization and tensor tubal ranks were introduced by Kilmer and her collaborators in \cite{KBHH13, KM11, KMP08, ZEAHK14}.  They are now widely used in engineering \cite{CXZ20, MQW20, MQW21, LYQX20, SHKM14, SNZ21, XCGZ21, XCGZ21a, YHHH16, ZSKA18, ZA17, ZHW21, ZLLZ18}.   In particular, Kilmer and Martin \cite{KM11} proposed T-SVD factorization.    By going through the Discrete Fourier Transform (DFT), matrix SVD and inverse DFT, a third order tensor is diagonalized to an f-diagonal third order tensor.  The tensor tubal rank is defined based upon such an f-diagonal tensor.     The matrix SVD should follow the standard non-increasing ordering for the singular values of the matrices involved.   If a different ordering is used, the diagonalization result would be different.

We call the above particular diagonalization  the Kilmer-Martin mapping, and say that two third order tensors are orthogonally equivalent if one of them can be obtained by the product of another with some orthogonal tensors.  We show that if two third order tensors are orthogonally equivalent, then their Kilmer-Martin mappings are the same.  Thus, the real f-diagonal tensor obtained by the Kilmer-Martin mapping of a third order tensor extracts the main features of that third order tensor.   We call the absolute values of the diagonal entries of the f-diagonal tensor as the singular values of the original third order tensor, and the number of the nonzero singular values as the T-rank of the third order tensor.    Some properties of singular values and T-ranks are studied.

The largest singular value of a real matrix is always greater than or equal to the absolute value of any entry of that matrix. We show that this is also true for third order tensors. Furthermore, based on this property,
we show that the best T-rank one approximation of a third order tensor can be given by its largest singular value and related orthogonal tensors.

The remaining of this paper is distributed as follows.  In the next section, some preliminary knowledge on T-product of third order tensors is reviewed.    The Kilmer-Martin mapping and orthogonal equivalence are defined in Section 3.  We show there that the Kilmer-Martin mappings of two orthogonally equivalent tensors are the same.   Singular values and T-ranks are defined in Section 4.   Their properties are also studied there.    In Section 5, we study the best T-rank one approximation of a third order tensor. 
Some further discussion is made in Section 6.

\section{Preliminaries}

In this paper, real matrices are denoted by capital Roman letters $A,B,\cdots$, complex matrices are denoted by capital Greek letters $\Delta, \Sigma, \cdots$, and tensors are denoted by Euler script letters $\mathcal{ A},\mathcal{B},\cdots$.    We use $\R$ to denote the real number field, and $\CC$ to denote the complex number field. For a third order tensor $\mathcal{A}\in \mathbb{R}^{m\times n \times p}$, its $(i,j,k)$-th element is represented by $a_{ijk}$, and use the Matlab notation $\mathcal{A}(i,:,:)$, $\mathcal{A}(:,i,:)$ and $\mathcal{A}(:,: ,i)$ respectively to represent the $i$-th horizontal, lateral and frontal slice of the $\mathcal{A}$. The frontal slice $\mathcal{A}(:,:,i)$ is represented by $A^{(i)}$. Define $\|\mathcal{A}\|_F:=
\sqrt{\sum_{ijk}|a_{ijk }|^2}$. 

For a third order tensor $\A \in \R^{m \times n \times p}$, 
as in  \cite{KBHH13,KM11}, define
$${\rm bcirc}(\A):= \left(\begin{aligned} A^{(1)}\ & A^{(p)} & A^{(p-1)} & \cdots & A^{(2)}\ \\ A^{(2)} & A^{(1)} & A^{(p)} & \cdots & A^{(3)}\\ \cdot\ \ \ & \ \cdot & \cdot\ \  & \cdots & \cdot\ \ \ \\
	\cdot\ \ \ & \ \cdot & \cdot\ \  & \cdots & \cdot\ \ \ \\
	A^{(p)} & A^{(p-1)} & A^{(p-2)} & \cdots & A^{(1)} \end{aligned}\right),$$
and bcirc$^{-1}($bcirc$(\A)):= \A$.

For a third order tensor $\A \in \R^{m \times n \times p}$, its transpose
is defined as
$$\A^\top = {\rm bcirc}^{-1}[({\rm birc}(\A))^\top].$$
This will be the same as the definition in \cite{KBHH13, KM11}.  The identity tensor $\I_{nnp}$ may also be defined as
$$\I_{nnp} = {\rm bcirc}^{-1}(I_{np}),$$
where $I_{np}$ is the identity matrix in $\R^{np \times np}$.  

A third order tensor $\S$ in $\R^{m \times n \times p}$ is f-diagonal in the sense of \cite{KBHH13, KM11} if all of its frontal slices $S^{(1)}, \cdots, S^{(p)}$ are diagonal.   We call the diagonal entries of  $S^{(1)}, \cdots, S^{(p)}$ as diagonal entries of $\S$.

For a third order tensor $\A \in \R^{m \times n \times p}$, it is defined \cite{KM11} that
$${\rm unfold}(\A) := \left(\begin{aligned} A^{(1)}\\ A^{(2)}\\ \cdot\ \ \\ \cdot\ \ \\ \cdot\ \ \\ A^{(p)}\end{aligned}\right) \in \R^{mp \times n},$$
and fold$($unfold$(\A)) := \A$.   For $\A \in \R^{m \times s \times p}$ and $\B \in \R^{s \times n \times p}$, the T-product of $\A$ and $\B$ is defined as
$\A * \B :=$ fold$(${bcirc$(\A)$unfold$(\B)) \in \R^{m \times n \times p}$.   Then, we see that
	\begin{equation} \label{e2.1}
		\A * \B = {\rm bcirc}^{-1}({\rm bcirc}(\A){\rm bcirc}(\B)).
	\end{equation}
	Thus, the bcirc and bcirc$^{-1}$ operations not only form a one-to-one relationship between third order tensors and block circulant matrices, but also their product operation is reserved.   By \cite{KM11}, the T-product operation (\ref{e2.1}) can be done by applying the fast Fourier transform (FFT).   The computational cost for this is $O(mnsp)$ flops.
	
	A tensor $\A \in \R^{n \times n \times p}$ has an inverse $\A^{-1} := \B \in \R^{n \times n \times p}$ if
	$$\A * \B = \B * \A = \I_{nnp}.$$
	If $\Q^{-1} = \Q^\top$ for $\Q \in \R^{n \times n \times p}$, then $\Q$ is called an orthogonal tensor.   

	\begin{definition} \label{d2.1}
		Suppose that $\A \in \R^{m \times n \times p}$.  The smallest integer $r$ such that
		\begin{equation} \label{e2.2}
			\A = \B * \C,
		\end{equation}
		where $\B \in \R^{m \times r \times p}$ and $\C \in \R^{r \times n \times p}$, is called the tensor tubal rank of $\A$.
	\end{definition}
	
	This definition was implicitly raised by Kilmer and Martin \cite{KM11} in 2011. In \cite{QY21}, this definition was formally used.
	
	\section{The Kilmer-Martin Mapping and Orthogonal Equivalence}
	\
	\newline
	Suppose that $\A \in \R^{m \times n \times p}$.  By (3.1) of \cite{KM11}, we may block-diagonalize bcirc$(\A)$ as
	\begin{equation} \label{e3.3}
		\Delta(\A) := (F_p \otimes I_m) {\rm bcirc}(\A) (F_p^* \otimes I_n) = \begin{bmatrix}
			\Delta^{(1)} &  &  & \\
			& \Delta^{(2)} &  & \\
			&  & \ddots& \\
			& &  & \Delta^{(p)} \\
		\end{bmatrix},
	\end{equation}
	where $F_p$ is the $p \times p$ DFT matrix,
	$F_p^*$ is its conjugate transpose, $\otimes$ denotes the Kronecker product, $\Delta^{(k)} \in \CC^{m \times n}$ for $k = 1, \cdots, p$.   For each matrix $\Delta^{(k)}$, compute its SVD
	$$\Delta^{(k)} = \Phi^{(k)}\Sigma^{(k)}{\Psi^{(k)}}^{*},$$
	where $\Phi^{(k)} \in \CC^{m \times m}$ and $\Psi^{(k)} \in \CC^{n \times n}$ are unitary matrices, $\Sigma^{(k)} \in \CC^{m \times n}$ is a diagonal matrix, the singular values of $\Delta^{(k)}$ follow the standard non-increasing order.  Denote
	\begin{equation} \label{e3.5}
		\Sigma(\A) :=  \begin{bmatrix}
			\Sigma^{(1)} &  &  & \\
			& \Sigma^{(2)} &  & \\
			&  & \ddots& \\
			& &  & \Sigma^{(p)} \\
		\end{bmatrix}.
	\end{equation}
	Let
	\begin{equation} \label{e3.6}
		\S = \S(\A) :=  {\rm bcirc}^{-1}\left((F_p^* \otimes I_m)\Sigma(\A)(F_p \otimes I_n)\right).
	\end{equation}
	Then $\S(\A) \in \R^{m \times n \times p}$ is a real f-diagonal tensor.  We call $\S(\cdot)$ the Kilmer-Martin mapping.  In particular, by Eq. (\ref{e3.6}), we have
	\begin{equation} \label{explicitform}
		\S(i, i, k) = {1 \over p}\sum_{l=1}^p \bar \omega^{(k-1)(l-1)}\Sigma(i, i, l),
	\end{equation}
	for $i = 1, \cdots, \min \{ m, n \}$, $k = 1, \cdots, p$, and  $\omega = e^{-2\pi\sqrt{-1}/p}$.
	Note that $\A$ has $mnp$ entries, and $\S = \S(\A)$ has $p\min \{ m, n\}$ diagonal entries.  In a certain sense, the main features of $\A$ are extracted in the diagonal entries of $\S$.
	
	As noticed in \cite{KM11}, the particular diagonalization $\S(\A)$ was achieved using the standard non-increasing ordering for the singular values of each $\Delta^{(k)}$.   If a different ordering is used, a different diagonalization $\S_1(\A)$ would be achieved.  Then the set of the diagonal entries of $\S_1(\A)$ can be different from the set of the diagonal entries of $\S(\A)$.
	
	Let
	$$\Phi(\A) :=  \begin{bmatrix}
		\Phi^{(1)} &  &  & \\
		& \Phi^{(2)} &  & \\
		&  & \ddots& \\
		& &  & \Phi^{(p)} \\
	\end{bmatrix},$$
	$$\Psi(\A) :=  \begin{bmatrix}
		\Psi^{(1)} &  &  & \\
		& \Psi^{(2)} &  & \\
		&  & \ddots& \\
		& &  & \Psi^{(p)} \\
	\end{bmatrix},$$
	$$\U = \U(\A) = {\rm bcirc}^{-1}\left((F_p^* \otimes I_m)\Phi(\A)(F_p \otimes I_n)\right),$$
	$$\V = \V(\A) = {\rm bcirc}^{-1}\left((F_p^* \otimes I_m)\Psi(\A)(F_p \otimes I_n)\right).$$
	Then $\U \in \R^{m \times m \times p}$ and $\V \in \R^{n \times n \times p}$ are orthogonal tensors, and $\A$ has its T-SVD
	\begin{equation} \label{e3.7}
		\A = \U * \S * \V^\top.
	\end{equation}
	
	Theorem 4.3 of \cite{KM11} showed that an Eckart-Young like theorem holds for the tensor tubal rank of $\A$ here.   For the Kilmer-Martin T-SVD factorization (\ref{e3.7}), by \cite{KM11, KMP08}, we have
	\begin{equation} \label{e3.8}
		\sum_{k=1}^p \S(1, 1, k)^2 \ge \sum_{k=1}^p \S(2, 2, k)^2 \ge \cdots  \ge  \sum_{k=1}^p \S(\min \{ m,n \}, \min \{ m,n \}, k)^2.
	\end{equation}
	
	Recently, Qi and Yu \cite{QY21} defined the $i$th largest T-singular value of $\A$ as
	$$\lambda_i := \sqrt{\sum_{k=1}^p \S(i, i, k)^2},$$
	for $i = 1, \cdots, \min \{ m, n \}$, and used T-singular values to define the tail energy for the error estimate of a proposed tensor sketching algorithm.    T-singular values are nonnegative numbers.   The number of the nonzero T-singular values of $\A$ is equal to the tensor tubal rank of $\A$.

	\begin{definition}
		Suppose that $\A, \B \in \R^{m \times n \times p}$. If there are orthogonal tensors $\Y \in \R^{m \times m \times p}$ and $\Z \in \R^{n \times n \times p}$ such that
		$$\A = \Y * \B * \Z^\top.$$
		Then we say that $\A$ and $\B$ are orthogonally equivalent.
	\end{definition}
	
	\begin{theorem} \label{t3.2}
		Suppose that $\A, \B \in \R^{m \times n \times p}$ are orthogonally equivalent,  $\A = \Y * \B * \Z^\top$, where $\Y \in \R^{m \times m \times p}$ and $\Z \in \R^{n \times n \times p}$ are orthogonal tensors.
		Then
		\begin{equation} \label{e3.9}
			\S(\A) = \S(\B).
		\end{equation}
	\end{theorem}
	{\bf Proof}
		We have
		$${\rm bicrc}(\A) = {\rm bcirc}(\Y){\rm bcirc}(\B){\rm bcirc}(\Z^\top).$$
		Apply $(F_p \otimes I_m)$ to the left and $(F_p^* \otimes I_n)$ to the right of each of the block circulant matrices in the above expression, where $F_p$ is the $p \times p$ DFT matrix, $F^*_p$ is its conjugate transpose, $\otimes$ denotes the Kronecker product.
		For each diagonal block, by SVD of matrices, we have
		$$\Delta(\A)^{(k)} = \Xi^{(k)} \Delta(\B)^{(k)}(\Theta^{(k)})^{*},$$
		where $\Xi^{(k)} \in \CC^{m \times m}$ and $\Theta^{(k)} \in \CC^{n \times n}$ are unitary matrices for $k = 1, \cdots, p$.  Then $\Delta(\A)^{(k)}$ and $\Delta(\B)^{(k)}$ have the same set of singular values for $k = 1, \cdots, p$.    This implies that
		$$\Sigma(\A)^{(k)} = \Sigma(\B)^{(k)},$$
		for $k = 1, \cdots, p$, i.e.,
		$$\Sigma(\A) = \Sigma(\B),$$
		which implies (\ref{e3.9}).
	
	\begin{corollary} \label{c3.3}
		Suppose that $\A, \B \in \R^{m \times n \times p}$ are orthogonally equivalent.   Then they have the same tensor tubal rank and T-singular value set.
	\end{corollary}
	
	Suppose that $\A \in \R^{m \times n \times p}$ is f-diagonal.  It is possible that $\S(\A)$ is very different from $\A$.  For example, let $\A = (a_{ijk}) \in \R^{3 \times 3 \times 3}$ be f-diagonal with  $a_{221} = 6$, $a_{112} = 5$, $a_{332} = 9$, $a_{333}=9$.  The other entries of $\A$ are zero.   Let $\S = \S(\A)$.   Then we have
	$$\S(1,1,1) = 12, \ \S(2,2,1) = 6, \  \S(3,3,1) ={ 5}, \S(1, 1, 2) = \S(1, 1, 3) = 3.$$
	The other entries of $\S$ are zero.    We see that $\S$ and $\A$ are very different.
	
	In \cite[Theorem 2.1]{LHPQ21}, the following result was proved.
	
	\begin{proposition} \label{p3.4}
		Suppose that $\A, \B \in \R^{m \times n \times p}$ are orthogonally equivalent,  $\A = \Y * \B * \Z^\top$, where $\Y \in \R^{m \times m \times p}$ and $\Z \in \R^{n \times n \times p}$ are orthogonal tensors.
		If $\B = (b_{ijk})$ is f-diagonal, then
		$$\max \{ |b_{ijk}| : 1 \le i \le m, 1 \le j \le n, 1 \le k \le p \} \le \sigma_1(\A),$$ where $\sigma_1(\A)$ is defined in Definition \ref{d4.1}.
	\end{proposition}
	
	By Theorem \ref{t3.2} and Proposition \ref{p3.4}, we have the following proposition.
	
	\begin{proposition} \label{p3.5}
		Suppose that $\A = (a_{ijk}) \in \R^{m \times n \times p}$ is f-diagonal.     Then
		$$\max \{ |a_{ijk}| : 1 \le i \le m, 1 \le j \le n, 1 \le k \le p \} \le \sigma_1(\A),$$ where $\sigma_1(\A)$ is defined in Definition \ref{d4.1}.
	\end{proposition}

	\section{Singular Values and T-Rank}
	
	\begin{definition} \label{d4.1}
		Suppose that $\A \in \R^{m \times n \times p}$.    Let (\ref{e3.7}) be the Kilmer-Martin T-SVD factorization of $\A$, where $\S = \S(\A)$.    The absolute values of the diagonal entries of the frontal slices of $\S$ are called the singular values of $\A$.  Let $s$ be a positive integer such that $1 \le s \le p\min \{ m, n \}$. The $s$th largest singular value of $\A$ is denoted as $\sigma_s(\A)$.    The number of the nonzero singular values of $\A$ is called the T-rank of $\A$.    Let
		$\S_s  \in \R^{m \times n \times p}$ be an f-diagonal tensor, such that its entries are the same as the entries of $\S$, where $\sigma_i(\A)$ for $1 \le i \le s$ are located, and its other entries are zero.   Denote $\A_s = \U * \S_s * \V^\top$.
	\end{definition}
	
	By Theorem \ref{t3.2} and Definition \ref{d4.1}, we have the following corollary.
	
	\begin{corollary}
		Suppose that $\A, \B \in \R^{m \times n \times p}$ are orthogonally equivalent.   Then they have the same T-rank and singular value set.
	\end{corollary}
	
	Zhang and Aeron \cite{ZA17} defined singular values of a third order tensor $\A \in \R^{m \times n \times p}$.   Suppose that $\A$ has a T-SVD
	$$\A = \U * \S * \V^\top,$$
	where $\U$ and $\V$ are orthogonal tensors, $\S$ is an f-diagonal tensor.   Then they call the entries of $\S$ the singular values of $\A$ \cite[Definition II.7]{ZA17}.   First, if not specifying $\S \not =  \S(\A)$, this definition is not well-defined.  Second, off-diagonal entries of $\S$ are zeros.  They are not needed to be involved.   Third, some diagonal entries of $\S$ may be negative.  Hence, our definition is different from theirs.
	
	\begin{proposition} \label{p4.3}
		Suppose that $\A = (a_{ijk}) \in \R^{m \times n \times p}$ has only one nonzero entry.   Then the T-rank of $\A$ is equal to one.
	\end{proposition}
{\bf Proof}
		Assume that $a_{i_0j_0k_0} = a \not = 0$, and the other entries of $\A$ are zero.  By writing out $F_p$ and $F^*_p$ explicitly, in (\ref{e3.3}), for $i = 1, \cdots, m$, $j = 1, \cdots, n$ and $k = 1, \cdots, p$, we have
		$$\Delta^{(k)}(i, j) = \sum_{l=1}^p \omega^{(k-1)(l-1)}A^{(l)}(i, j),$$
		where $\omega = e^{-2\pi\sqrt{-1}/p}$.
		Then for $k = 1, \cdots, p$,
		$$\Delta^{(k)}(i_0, j_0) = \omega^{(k-1)(k_0-1)}A^{(k_0)}(i_0, j_0).$$
		The other entries of $\Delta(\A)$ are zero.  Consider the SVD of $\Delta^{(k)}$.  The singular values of  $\Delta^{(k)}$ are the square roots of eigenvalues of $(\Delta^{(k)})^*\Delta^{(k)}$.  Then, $(\Delta^{(k)})^*\Delta^{(k)}$ only has a nonzero entry
		$$(\Delta^{(k)})^*\Delta^{(k)}(j_0,j_0) = a_{i_0j_0k_0}^2,$$
		for $k = 1, \cdots, p$.  Then, for $k= 1, \cdots, p$, $\Delta^{(k)}$ has only one nonzero singular value $|a|$.   This implies that $\Sigma^{(k)}(1, 1)  = |a_{i_0j_0k_0}| = |a|$ and $\Sigma^{(k)}(i, i) = 0$ for $k = 1, \cdots, p$ and $i = 2, \cdots, \min \{ m, n \}$.    By the inverse DFT, we have
		$$\S(1,1,1) = {1 \over p}\sum_{k=1}^p |a| = |a|,$$
		$$\S(1, 1, k) = {1 \over p}\sum_{l=1}^p \bar \omega^{(k-1)(l-1)} \Sigma^{(l)}(1, 1) = {|a| \over p}\sum_{l=1}^p \bar \omega^{(k-1)(l-1)} = 0,$$
		as $\sum_{l=1}^p \bar \omega^{(k-1)(l-1)} = 0$, for $k = 2, \cdots, p$.
		We also have
		$$\S(i, i, k) = 0$$
		for $i \ge 2$ and $k = 1, \cdots, p$.    Then $\S(\A)$ has only one nonzero entry $\S(1, 1, 1) = |a|$.

	\begin{proposition}  \label{p7.2}
		Suppose $\A \in \R^{m \times n \times p}$, $\A = \A' + \A''$.  Then
		$$\sigma_1(\A) \le \sigma_1(\A') + \sigma_1(\A'').$$
	\end{proposition}
{\bf Proof}
		Denote $A^{(i)}=\A(:,:,i)$, $A'^{(i)}=\A'(:,:,i)$ and $A''^{(i)}=\A''(:,:,i)$, then $A^{(i)}=A'^{(i)}+A''^{(i)}$. We have $$\mbox{bcirc}(\A)=\mbox{bcirc}(\A')+\mbox{bcirc}(\A'').$$
		Applying FFT to both sides, the above equation is transformed to the following:
		\[
		\begin{bmatrix}
			\Delta^{(1)} &  &  & \\
			& \Delta^{(2)} &  & \\
			&  & \ddots& \\
			& &  & \Delta^{(p)} \\
		\end{bmatrix}=
		\begin{bmatrix}
			\Delta'^{(1)} &  &  & \\
			& \Delta'^{(2)} &  & \\
			&  & \ddots& \\
			& &  & \Delta'^{(p)} \\
		\end{bmatrix}
		+\begin{bmatrix}
			\Delta''^{(1)} &  &  & \\
			& \Delta''^{(2)} &  & \\
			&  & \ddots& \\
			& &  & \Delta''^{(p)} \\
		\end{bmatrix},
		\]
		where $\Delta^{(i)}=\sum_{l=1}^{p}\omega^{(l-1)(i-1)}A^{(i)}=\sum_{l=1}^{p}\omega^{(l-1)(i-1)}(A'^{(i)}+A''^{(i)})=\Delta'^{(i)}+\Delta''^{(i)}$.
		
		Denote $\sigma_1(A)$ as the largest singular value of a matrix $A$. Then for each $i=1,\cdots,p$, we have
		$$\sigma_1(\Delta^{(i)})\le \sigma_1(\Delta'^{(i)}) +\sigma_1(\Delta''^{(i)}).$$
		Thus,
		\begin{eqnarray*}
			\max_{i=1,\cdots,p}\{\sigma_1(\Delta^{(i)})\}
			& \le & \max_{i=1,\cdots,p}\{\sigma_1(\Delta'^{(i)})\}+\max_{i=1,\cdots,p}\{\sigma_1(\Delta''^{(i)})\}
		\end{eqnarray*}
		Hence,
		\begin{eqnarray*}
			\sigma_1(\A) &\le& \sigma_1(\A')+\sigma_1(\A'').
		\end{eqnarray*}
	

	
	
	\section{The Best T-Rank One Approximation to a Third Order Tensor}
	In this section, we propose the best T-rank one approximation to a third order tensor. Before that,  we first give a theorem to show the relation between the largest singular value and the maximum absolute value of elements in a third order tensor.

	Recall that by (\ref{explicitform}),  for any  $\A = (a_{ijk}) \in \R^{m \times n \times p}$, we have
	\begin{equation} \label{n2.2}
		a_{ijk} = {1 \over p}\sum_{l=1}^p \bar \omega^{(k-1)(l-1)}\Delta(i, j, l),
	\end{equation}
	for $i = 1, \cdots, m$, $j=1, \cdots, n$, $k = 1, \cdots, p$.
	Furthermore, according to  (\ref{explicitform}), we have
	\begin{equation} \label{n2.3}
		\S(1, 1, 1) = {1 \over p}\sum_{l=1}^p \Sigma(1,1,l).
	\end{equation}
	Thus,  if not all $\Sigma(1,1,l)$ are equal, i.e., $\min_l \{ \Sigma(1,1,l) \} < \max_l \{ \Sigma(1,1,l) \}$, then
	$$\S(1,1,1) <  \max_l \{ \Sigma(1,1,l) \}.$$
	Hence, we have the following theorem.
	
	\begin{theorem}\label{th6.1}
		Suppose $\A = (a_{ijk}) \in \R^{m \times n \times p}$.   Then
		\begin{equation}
			\sigma_1(\A) \equiv \S(1, 1, 1) \ge \max \{ |a_{ijk}| : 1 \le i \le m, 1 \le j \le n, 1 \le k \le p \}.
		\end{equation}
	\end{theorem}
	
{\bf Proof}
		By the property of matrix singular values \cite[(2.3.8) and Corollary 2.4.3]{GV13}\footnote{\cite[(2.3.8) and Corollary 2.4.3]{GV13} is for real matrices, but they can be extended to complex matrices easily.}, we have
		\begin{equation}
			|\Delta(i,j,k)| \le \Sigma(1,1,k),
		\end{equation}
		for $i = 1, \cdots, m$, $j = 1, \cdots, n$, $k = 1, \cdots, p$.
		Then by (\ref{n2.2}) and (\ref{n2.3}), we have
		\begin{eqnarray*}
			|a_{ijk}| & \le & {1 \over p}\sum_{l=1}^p  | \bar \omega^{(k-1)(l-1)}\Delta(i, j, l)|\\
			& \le & {1 \over p}\sum_{l=1}^p  |\Delta(i, j, l)|\\
			& \le & {1 \over p}\sum_{l=1}^p \Sigma(1,1,l)\\
			& = & \S(1, 1, 1),
		\end{eqnarray*}
		for $i = 1, \cdots, m$, $j = 1, \cdots, n$, $k = 1, \cdots, p$, i.e,
		\begin{equation*}
			\sigma_1(\A) \equiv \S(1, 1, 1) \ge \max \{ |a_{ijk}| : 1 \le i \le m, 1 \le j \le n, 1 \le k \le p \}.
		\end{equation*}

	Such a property holds in the matrix case \cite{GV13}.   Our theorem shows that it holds as well in the third order tensor case. Based on Theorem \ref{th6.1}, we have the following theorem.
	
	\begin{theorem} \label{p6.1}
		For any $\A \in \R^{m \times n \times p}$, $\A_1$ is the best T-rank one approximation of $\A$, where ${\A}_1$ is defined by in Definition \ref{d4.1}.
	\end{theorem}
	
{\bf Proof}
		Let $\A \in \R^{m \times n \times p}$.
		We have
		$$\|\A - \A_1\|_F^2 =  \| \S - \S_1 \|_F^2 = \sum_{i=2}^{p\min \{ m, n \}} \sigma_i^2.$$
		Now, assume that $\B \in \R^{m \times n \times p}$ has T-rank one.  Suppose that $\B$ has a Kilmer-Martin T-SVD factorization
		$$\B = \Y * \D * \Z^\top,$$
		where $\Y \in \R^{m \times m \times p}$ and $\Z \in \R^{n \times n \times p}$ are orthogonal, $\D \in \R^{m \times n \times p}$ is f-diagonal and has at most $1$ nonzero elements.  Let
		$$\A' = \Y^\top * \A * \Z.$$
		Then
		$$\| \A - \B \|_F^2 = \| \A' - \D \|_F^2 \ge \|\A'\|_F^2 - \sigma_1(\A')^2.$$
		where the inequality follows from Theorem \ref{th6.1}  applied to $\A'$.
		
		By Theorem \ref{t3.2}, $\sigma_i(\A') = \sigma_i(\A)$ for $i = 1, \cdots, \min \{ m, n \}$.
		
		Then we have
		$$\| \A - \B \|_F^2 \ge \|\A'\|_F^2 -  \sigma_1(\A')^2 = \|\A\|_F^2 - \sigma_1(\A)^2 = \sum_{i=2}^{p\min \{ m, n \}} \sigma_i^2 = \|\A - \A_1\|_F^2.$$
		This shows that ${\A}_1$ is the best T-rank one approximation of $\A$.

	\section{Further Discussion}
	\begin{enumerate}
		\item Comparing with the tensor tubal rank, the T-rank is simple in the best rank one approximation to third order tensors.
		\item The T-rank is not sub-additive.   Let $\B = (b_{ijk}), \C = (c_{ijk}), \D = (d_{ijk}), \E = (e_{ijk}) \in \R^{3 \times 3 \times 3}$, and each of them has exactly one nonzero entry as $b_{221} = 6$, $c_{112} = 5$, $d_{332} = 9$, $e_{333}=9$.  Their other entries are zero.   Then by Proposition \ref{p4.3}, they are all T-rank one tensors.   Let $\A = \B + \C + \D + \E$.  Then we find T-rank$(\A) = 5 >$ T-rank$(\B)+$ T-rank$(\C)+$ T-rank$(\D)+$  T-rank$(\E) = 4$.   In fact, for $\S = \S(\A)$, we have
		$$\S(1,1,1) = 12, \ \S(2,2,1) = 6, \  \S(3,3,1) ={5}, \S(1, 1, 2) = \S(1, 1, 3) = 3.$$
		The other entries of $\S$ are zero.
		\item  Kilmer and Martin \cite[Theorem 4.3]{KM11} showed that an Eckart-Young like theorem holds for the tensor tubal rank of third order tensors.    Does another Eckart-Young like theorem hold for the T-rank of third order tensors?   This may be an interesting point for further exploration.
	\end{enumerate}
	\vspace{3mm}
	{\bf Acknowledgement}
	
	We are thankful to Mr. Xin Chen for his comments.

\bigskip


\end{document}